\documentclass[a4paper,12pt]{article}
\usepackage[colorlinks=true, linkcolor=blue, citecolor=blue]{hyperref}
\usepackage[margin=1.55cm]{geometry}
\usepackage{helvet}
\usepackage{amsmath}
\usepackage{amssymb}
\usepackage{mathtools} %loads amsmath as well
\usepackage{enumerate}% http://ctan.org/pkg/enumerate
\usepackage{graphicx}
\usepackage{float}
\usepackage[superscript, biblabel, nomove]{cite}
\usepackage{color}
\usepackage{url}
\usepackage{subcaption}
\usepackage{caption}
\usepackage{caption}
\usepackage{hyperref}
\hypersetup{pdfstartview={XYZ null null 1.00}}
\usepackage{amsmath}
\usepackage{amssymb}
\usepackage{mathtools} %loads amsmath as well
\usepackage{enumerate}% http://ctan.org/pkg/enumerate\jmath{{{{k}}}}
\usepackage{graphicx}
\usepackage{float}
\usepackage{color}
\usepackage{url}
\usepackage{subcaption}
\usepackage{caption}
\usepackage{hyperref}
\usepackage{amssymb}
\usepackage{amsmath}
\usepackage{amsthm}
\usepackage{float}
\usepackage{siunitx}
\usepackage{booktabs}
\usepackage[english]{babel}
\usepackage[nottoc]{tocbibind}
\usepackage{lipsum}
\usepackage{booktabs}
\usepackage{siunitx}
\usepackage{makecell}
\usepackage{changepage}
\usepackage{orcidlink}
\newcommand\blfootnote[1]{%
  \begingroup
  \renewcommand\thefootnote{}\footnote{#1}%
  \addtocounter{footnote}{-1}%
  \endgroup
}
\usepackage[linewidth=2pt]{mdframed}

\newmdenv[
topline=false,
bottomline=false,
rightline=false,
skipabove=\topsep,
skipbelow=\topsep,
linewidth=4
]{siderules}
\usepackage{hyperref}
\usepackage{xcolor}
\usepackage[most]{tcolorbox}
\usepackage[nottoc]{tocbibind}

\usepackage{hyperref}
\definecolor{DarkBlue}{rgb}{0.00,0.00,0.50}
\hypersetup{colorlinks=true,
            linkcolor=blue,
            filecolor=magenta,
            urlcolor=cyan,
            citecolor=DarkBlue,
            pdftitle={MSc UoR Computer Science Report},
            }

\usepackage{authblk}
\usepackage{blindtext}
\usepackage{datetime}
\newdateformat{monthdayyeardate}{%
  \monthname[\THEMONTH]~\THEDAY, \THEYEAR}
  \numberwithin{equation}{section}
\newtheorem{definition}{Definition}
\let\olddefinition\definition
\renewcommand{\definition}{\olddefinition\normalfont}
\newtheorem{prb}{Problem}
\let\oldprb\prb
\renewcommand{\prb}{\oldprb\normalfont}

\title{Iterative Procedure for Non-Linear Fractional  Integro-Differential Equations via Daftardar--Jafari Polynomials} 
\author {Qasim Khan${}^{\footnote{Corresponding author.}}$}
\author {Anthony Suen}
\affil{ Department of Mathematics and Information Technology,The Education University of Hong Kong,	10 Lo Ping Road, Tai Po, N.T,	Hong Kong. }
\date{}
\begin{document}
	\maketitle
\blfootnote{E-mail addresses: {\href{mailto:qasimkhan@s.eduhk.hk}{qasimkhan@s.eduhk.hk} (Q. Khan),    \href{mailto:acksuen@eduhk.hk}{acksuen@eduhk.hk}  (A. Suen). }} \blfootnote{\emph{Preprint submitted to Elsevier}}
\begin{abstract}
 In this paper, we introduce a novel approach called the Iterative Aboodh Transform Method (IATM) which utilizes Daftardar--Jafari polynomials for solving non-linear problems. Such method is employed to derive solutions for non-linear fractional partial integro-differential equations (FPIDEs). The key novelty of the suggested method is that it can be used for handling solutions of non-linear FPIDEs in a very simple and effective way. {More precisely, we show that Daftardar--Jafari polynomials have simple calculations as compared to Adomian polynomials with higher accuracy}. The results obtained within the Daftardar--Jafari polynomials are demonstrated with graphs and tables, and the IATM's absolute error confirms the higher accuracy of the suggested method.
\end{abstract}
\section{Introduction}
Fractional partial integro-differential equations (FPIDEs), which incorporate both integrals and fractional derivatives, play a crucial role in the study of fractional calculus. A significant aspect of { FPIDEs} is the application of fractional specifications, which has become one of the essential mathematical tools for modelling and analysing real-world problems in natural sciences, engineering, and technology \cite{l2,6q,7q,8q}. Due to the non-local property of fractional derivative operators, the study on FPIDEs is highly non-trivial and challenging. Nonetheless, various researchers have made significant strides in studying FPIDEs analytical and numerically, here we briefly recall some results from the literature. Hassan et al. \cite{h2} employed the Chebyshev Wavelet Method (CWM) to tackle higher-order fractional integro-differential equations (FIDEs), demonstrating the method's potential in handling complex fractional dynamics. Mohyud-Din et al. \cite{h1} applied the CWM to non-linear FIDEs, showcasing its versatility across different types of equation. To derive analytical solutions for FPIDEs, Hussain et al. \cite{22} utilized the variation iteration approach, while Mittal et al. \cite{23} implemented the Adomian decomposition method to solve FPIDEs. Additionally, Awawdeh et al. \cite{21} employed the homotopy analysis method to derive analytical solutions for linear FPIDEs. Eslahchi et al. explored the Jacobi method for solving non-linear FPIDEs, investigating stability and convergence in their analysis \cite{25}. Zhao et al. \cite{26} introduced a piecewise polynomial collocation approach to address FPIDEs with weakly singular kernels, illustrating an innovative numerical strategy. Rawashdeh proposed a numerical method leveraging polynomial splines for effective numerical solutions of FPIDEs \cite{24}. Furthermore, Unhale et al. \cite{28} suggested a collocation method using shifted Legendre polynomials and Chebyshev polynomials to solve non-linear FPIDEs, while Avazzadeh et al. combined Legendre wavelets with a fractional integration operational matrix in a hybrid approach \cite{29}. Lastly, Arshed applied the B-spline method to find solutions to linear FPIDEs \cite{27}, and Dehestani et al. introduced the {collocation} method for the numerical simulation of FPIDEs \cite{30}. This evolving body of research underscores the growing significance of FPIDEs and the diverse methodologies being developed to address them effectively.

{On the other hand, It is also worth mentioning a special type of fractional equations, which are known as Stochastic fractional integro-differential equations (SFIDEs). SFIDEs are mathematical models that integrate the principles of stochastic processes, fractional calculus, and integral equations, enabling the representation of systems with inherent randomness and memory effects. SFIDEs are particularly useful in fields such as finance, physics, and biology, where they model phenomena characterized by unpredictable behaviour and long-term dependencies, such as stock price movements or the diffusion of particles in complex media. References and related discussions can be found in Mirzaee et al.\cite{mirzaee2024meshless, mirzaee2021approximate, mirzaee2021bicubic, mirzaee2020cubic}, Solhi et al. \cite{solhi2024enhanced}, Alipour and Mirzaee \cite{alipour2020iterative}, Mirzaee and Alipour \cite{mirzaee2019numerical}, Samadyar and Mirzaee\cite{samadyar2019numerical}.}

In the current study, we introduce the Iterative Aboodh transform method (IATM) to obtain the approximate solution of non-linear FPIDEs. To represent fractional derivatives, the Caputo derivative operator is utilised and the Daftardar--Jafari polynomials from \cite{l5,qa} are used to express the non-linear terms in each targeted problem. The Daftardar--Jafari polynomials have a straightforward implementation in a given infinite series form that can provide the sufficient degree of accuracy. The solutions of three examples are then illustrated in details to confirm the reliability of the suggested technique. Using MAPLE software, a very simple algorithm of the present method is constructed with the aid of Daftardar--Jafari polynomials for non-linear FPIDEs.

The IATM is comparably simple in the sense that it involves fewer calculations, which makes it suitable for extending to the study on solutions of other FPIDEs and their associated systems. { More precisely, this
 research work has the following multiple advantages:
\begin{itemize}
  \item We successfully achieve highly accuracy with the same number of iterations and parameters than the results recently published in Khan et al. \cite{LADM} and the related works cited therein.
  \item The non-linear terms are controlled by Daftardar--Jafari polynomials, which allow us to incorporate non-linearity directly into a more simple algebraic expression without involving derivatives as compare with Adomian polynomials and He's polynomials.
  \item  The mathematical model investigated has significant physical interpretation on various subjects, such as the heat flow in materials with memory and the mechanics of viscoelasticity. Reference \cite{ex1} and the references cited therein will provide further insights into these applications, which can illustrate the model's utility for understanding complex material behaviours.
    \item  We deal with one of the most commonly used fractional derivative operators, known as the Caputo derivative, which can improve the modelling accuracy for describing physical phenomena with viscoelastic forces; refer to equation \eqref{s33}.
    \item We believe that the newly proposed IATM may bring insights on the study of other related fractional equations, which include SFIDEs as discussed before. Further investigations will be carried out in subsequent research projects later.
\end{itemize}
}

{ The rest of the paper is organised as follows. In Section \ref{sec1}, we give the definitions of some fundamental concepts related to fractional calculus and operators. In Section \ref{sec3}, we introduce the IATM for solving non-linear FPIDEs. In Section \ref{sec4}, we provide the details on the numerical implementation of the proposed method, with the showcases of some specific problems and their solutions. In Section \ref{sec5}, we discuss the results of our numerical experiments and present some further implications. Finally, we conclude in Section \ref{sec6} with a summary of findings and reflections on the applicability of the proposed approach.}
\section{Basic Definitions}\label{sec1}
{In this section, we provide some definitions that are essential for implementing the iterative scheme. Interested readers are encouraged to consult the references \cite{d1,d2,ML,AT,AT1,ATF} for further information. To begin with, we define the Caputo operator for fractional derivatives.}
\begin{definition}
The Caputo operator for fractional derivatives of order ${{{{\alpha}}}}$ is given by
\begin{equation}\label{caputo eqn}
D_{{{{{t}}}}}^{{{{{\alpha}}}}}{{{{u}}}}({{{{x}}}},{{{{t}}}})=
\frac{1}{\Gamma({{{j}}}-{{{{\alpha}}}})}\int^{{{{t}}}}_0({{{{t}}}}
-{{{{t}}}}_{o})^{{{{j}}}-{{{{\alpha}}}}-1}\frac{\partial^{{{j}}} {{{{u}}}}({{{{x}}}},{{{{t}}}})}{\partial {{{{t}}}}^{{{j}}}}d{{{{t}}}}, ~~~~ if~~~ {{{j}}}-1<{{{{\alpha}}}} \leq {{{j}}}.
\end{equation}
If ${{{{\alpha}}}}={{{j}}}$ for ${{{j}}} \in \mathbb{N}$, then the operator $D_{{{{{t}}}}}^{{{{{\alpha}}}}}$ is reduced to the ordinary time derivative
\begin{equation}\nonumber
D_{{{{{t}}}}}^{{{{{\alpha}}}}}{{{{u}}}}({{{{x}}}},{{{{t}}}})=\frac{\partial^{{{j}}} {{{{u}}}}({{{{x}}}},{{{{t}}}})}{\partial {{{{t}}}}^{{{j}}}}.
\end{equation}
  \end{definition}
  Next, we introduce the well-known Laplace transform $\mathcal{L}(\cdot)$ and its interactions with Caputo operator.
  \begin{definition}
 For a function $g({{{{{t}}}}})$, the Laplace transform is given by
\[
{G(s)}=\mathcal{L}[{g({{{{t}}}})}]=\int_{0}^{\infty}e^{-s{{{t}}}}{g}({{{{{t}}}}})d{{{{t}}}}.
\]
We further have the following identities involving Laplace transform and Caputo operator:
\begin{equation*}
\begin{split}
    \mathcal{L} (D_{{{{t}}}}^{{{{{\alpha}}}}}{{{{u}}}}({{{{t}}}}))&=s^{{{{\alpha}}}} \mathcal{L}[{{{{u}}}}({{{{t}}}})]-\sum_{k=0}^{{j}-1}s^{{{{{{\alpha}}}}}-j-1}U^k(0^+), \ \ {j}-1<{{{{\alpha}}}}\leq {j}, \ \ \forall {j} \in \mathbb{N}. \\
 \mathcal{L} (D_{{{{t}}}}^{{{{{\alpha}}}}}{{{{u}}}}({{{{t}}}}))&=s^{{{{\alpha}}}} U(s)-s^{{{{{\alpha}}}}-1}U(0), \ \ \ 0<{{{{\alpha}}}}\leq {1} \ \ \ \ \normalfont{for} \ \ {j}=1.
\end{split}
\end{equation*}
 \end{definition}
{ We now introduce the Daftardar--Jafari polynomials (also known as the Daftardar--Jafari method) which is used for expressing the non-linear term in an approximate problem. Details are given as follows.
  \begin{definition}
  Suppose that $u(x,t)$ is a solution to the following functional equation
  \begin{equation*}
  u(x,t)=N(u(x,t))+f(x,t),
  \end{equation*}
  where $N(\cdot)$ is a non-linear operator between Banach spaces $\mathcal{B}_1$ and $\mathcal{B}_2$, and $f$ is a given function. We are looking for $u(x,t)$ in a series form
  \begin{equation*}
  u(x,t)=\sum_{j=0}^\infty u_j(x,t).
  \end{equation*}
The non-linear operator $N(\cdot)$ can then be decomposed as follows: 
\begin{equation}\label{s4}
N \left(\sum_{{j}=0}^\infty{{{{{u}}}}_{j}({{{{x}}}},{{{{t}}}})}\right)=N({{{{{u}}}}_0({{{{x}}}},{{{{t}}}})})+\sum_{{j}=0}^\infty \bigg[{N} \bigg(\sum_{i=0}^{{j}} {{{{u}}}}_i({{{{x}}}},{{{{t}}}})\bigg)-{N} \bigg(\sum_{i=0}^{{j}-1 }{{{{u}}}}_i({{{{x}}}},{{{{t}}}})\bigg)\bigg],
\end{equation}
where $u_0:=f$ and $u_1:=N(u_0)$ are respectively the Daftardar--Jafari polynomials $u_j$ when $j=0$ and $\j=1$. For $j\ge2$, the Daftardar--Jafari polynomials $u_j$ with respect to the solution $u$ are given by
\begin{equation*}
u_{j}:=N(u_0+\cdots+u_{j-1})-N(u_0+\cdots+u_{j-2}).
\end{equation*}
Throughout this paper, the Daftardar--Jafari polynomials are utilized to represent the non-linear terms in various problems, which will be explained in Section~\ref{sec3} and Section~\ref{sec4} later.
\end{definition}}
Lastly, we introduce the Aboodh transform (AT) and some related identities.
  \begin{definition}
Given $M\in(0,\infty)$ and $k_1,k_2\in(0,\infty]$, we define a set of functions {$\mathbf{S}:=\mathbf{S}(M,k_1,k_2)$} as follows
\begin{equation*}
{\mathbf{S}}=\{{{{{u}}}}({{{{t}}}}): |{{{{u}}}}({{{{t}}}})| < Me^{-st},\ \ t\ge0,\ \ k_1\le s\le k_2\}.
\end{equation*}
For any $u\in\mathbf{S}$, the Aboodh transform (AT) $\mathcal{A}[u]$ of $u$ is given by
\begin{equation*}
{k}({s}):=\mathcal{A}\left[{{{{{u}}}}}\right]=\frac{1}{s}\int_{0}^{\infty}{{{{{u}}}}}({{{{{t}}}}})e^{-s{{{{{t}}}}}}{d}{{{{{t}}}}},\ \ {{{{{t}}}}}\geq0,\ \ {{{k}}}_{1}\leq{s}\leq{{{k}}}_{2}.
\end{equation*}
The operator $\mathcal{A}$ enjoys the following properties:
\begin{itemize}
\item $\mathcal{A}\left[{{{{{{u}}}}}^{'}}({{{{{t}}}}})\right]=s\mathcal{A}[u](s)-\frac{{{{{u}}}}(0)}{s}, $
\item $\mathcal{A}\left[{{{{{{u}}}}}^{''}}({{{{{t}}}}})\right]=s^{2}\mathcal{A}[u](s)-\frac{{{{{{{u}}}}}^{'}}(0)}{s}-{{{{{u}}}}(0)}, $\\
\vdots
\item $\mathcal{A}\left[{{{{{{u}}}}}^{(n)}}({{{{{t}}}}})\right]=s^{n}\mathcal{A}[u](s)-\sum_{{j}=0}^{n-1} \frac{{{{{{{u}}}}}^{(n)}}(0)}{s^2-n+j}$.
\end{itemize}
Furthermore, if {$a$ and $b$ are constants and ${{{{u}}}}({{{{t}}}}),{{{{v}}}}({{{{t}}}})\in\mathbf{S}$, then we have}
\begin{equation*}
  \mathcal{A}\left[a{{{{{{u}}}}}}({{{{{t}}}}}) \pm b{{{{v}}}}({{{{t}}}}) \right]=a\mathcal{A}\left[{{{{{{u}}}}}}({{{{{t}}}}}) \right] \pm b\mathcal{A}\left[{{{{v}}}}({{{{t}}}})\right].
\end{equation*}

We also define $\mathcal{A}^{-1}$ as follows. If $f(t)$ is piecewise continuous and of exponential order for $t\ge0$ such that $\mathcal{A}[f(t)] = F(s)$, then $f(t)$ is called inverse Aboodh transform of $F(s)$ and we write $\mathcal{A}^{-1}[f(s)]=F(t)$.
 \end{definition}
{
One of the important applications of AT is for converting Caputo fractional differential operators into algebraic equations. Specifically, for $\alpha \in (0,1]$, the AT of $D^{\alpha}_{{{{{{t}}}}}}{{{{u}}}}$ is defined as
$$\mathcal{A}[D^{\alpha}_{{{{{{t}}}}}}{{{{u}}}}] ({{{{{x}}}}}, {{{{{s}}}}} ) = {{s}}^{\alpha}\mathcal{A}[{{{{{u}}}}}]({{{{{x}}}}},{{s}})-{{s}}^{\alpha-2}\mathcal{A}[{{{{u}}}}]({{{{{x}}}}},0).$$}

\section{IATM for  non-linear FPIDEs}\label{sec3}
In this section, we introduce the new approach known as the Iterative Aboodh Transform Method (IATM). In order to understand IATM, we consider the following non-linear FPIDEs 
\begin{equation}\label{s33}
\begin{split}
{{D}^{{{{\alpha}}}}_{{{{{t}}}}}{{{{{u}}}}({{{{x}}}},{{{{t}}}})}}+{{{{{u}}}}({{{{x}}}},{{{{t}}}})}{\frac {\partial}{{\partial}{{{{{x}}}}}}}{{{{{u}}}}({{{{x}}}},{{{{t}}}})}=\int_{0}^{{{{{t}}}}}\!
   \left( {{{{{t}}}}}-{q} \right) ^{{{{{\alpha}}}}-1}\,{\frac {\partial^2}{{\partial} {{{{{x}}}}}^2 }}
{{{{{u}}}}({{{{x}}}},{{{{t}}}})} {\rm d}q+f({{{{x}}}},{{{{t}}}})
   \ \ \ {{{{\alpha}}}}\in (0,1]
\end{split}
\end{equation}
where $f({{{{x}}}},{{{{t}}}})$ is a given source term and ${D}^{{{{\alpha}}}}_{{{{{t}}}}}$ is the Caputo type  fractional order derivative as defined in \eqref{caputo eqn}. The above equation \eqref{s33} is equipped with initial and boundary conditions
\begin{equation*}
   {{{{u}}}}({{{{x}}}},0)=g({{{{x}}}}),\ \
  {{{{u}}}}(0,{{{{t}}}})={{{{u}}}}(L,{{{{t}}}})=0.
\end{equation*}
Applying Aboodh transform (AT) on equation \eqref{s33}, we get
\begin{equation*}
\begin{split}
\mathcal{A}[{D}^{{{{\alpha}}}}_{{{{{t}}}}}{{{{u}}}}({{{{x}}}},{{{{t}}}})]=\mathcal{A}\left[\int_{0}^{{{{{t}}}}}\!
    \left( {{{{{t}}}}}-{q} \right) ^{{{{{\alpha}}}}-1}\,{\frac {\partial^2}{{\partial} {{{{x}}}}^2 }}
{{{{{u}}}}({{{{x}}}},{{{{t}}}})}{\rm d}q-{{{{{u}}}}({{{{x}}}},{{{{t}}}})}{\frac {\partial}{{\partial}{{{{x}}}}}}{{{{{u}}}}({{{{x}}}},{{{{t}}}})}+f({{{{x}}}},{{{{t}}}})\right].
\end{split}
\end{equation*}
Using the algebraic property of AT on fractional derivatives, equation \eqref{s33} can be rewritten as follows
\begin{equation}\label{oi}
{s^{{{{\alpha}}}}}\mathcal{A}\left[{{{{{{u}}}}}({{{{x}}}},{{{{t}}}})}\right]=s^{{{{{\alpha}}}}-2}{{{{{u}}}}}({{{{x}}}},0)+\mathcal{A}\left[\int_{0}^{{{{{t}}}}}\!
    \left( {{{{{t}}}}}-{q} \right) ^{{{{{\alpha}}}}-1}\,{\frac {\partial^2}{{\partial} {{{{x}}}}^2 }}
{{{{{u}}}}({{{{x}}}},{{{{t}}}})}{\rm d}q-{{{{{u}}}}({{{{x}}}},{{{{t}}}})}{\frac {\partial}{{\partial}{{{{x}}}}}}{{{{{u}}}}({{{{x}}}},{{{{t}}}})}+f({{{{x}}}},{{{{t}}}})\right],
\end{equation}
and after simplifying equation \eqref{oi}, we further obtain
\begin{equation}\label{s41}
\begin{split}
\mathcal{A}\left[{{{{{{u}}}}}({{{{x}}}},{{{{t}}}})}\right]=s^{-2}{{{{{u}}}}}({{{{x}}}},0)+\frac{1}{s^{{{{\alpha}}}}}\mathcal{A}\left[\int_{0}^{{{{{t}}}}}\!
    \left( {{{{{t}}}}}-{q} \right) ^{{{{{\alpha}}}}-1}\,{\frac {\partial^2}{{\partial} {{{{x}}}}^2 }}
{{{{{u}}}}({{{{x}}}},{{{{t}}}})}{\rm d}q-{{{{{u}}}}({{{{x}}}},{{{{t}}}})}{\frac {\partial}{{\partial}{{{{x}}}}}}{{{{{u}}}}({{{{x}}}},{{{{t}}}})}+f({{{{x}}}},{{{{t}}}})\right].
\end{split}
\end{equation}
We apply the inverse AT $\mathcal{A}^{-1}$ to equation \eqref{s41} and give
\begin{equation}\label{s51}
\begin{split}
{{{{{{u}}}}}({{{{x}}}},{{{{t}}}})}={{{{{u}}}}}({{{{x}}}},0)+\mathcal{A}^{-1}\left[\frac{1}{s^{{{{\alpha}}}}}\mathcal{A}\left[\int_{0}^{{{{{t}}}}}\!
    \left( {{{{{t}}}}}-{q} \right) ^{{{{{\alpha}}}}-1}\,{\frac {\partial^2}{{\partial} {{{{x}}}}^2 }}
{{{{{u}}}}({{{{x}}}},{{{{t}}}})}{\rm d}q-{{{{{u}}}}({{{{x}}}},{{{{t}}}})}{\frac {\partial}{{\partial}{{{{x}}}}}}{{{{{u}}}}({{{{x}}}},{{{{t}}}})}+f({{{{x}}}},{{{{t}}}})\right]\right].
\end{split}
\end{equation}
The assumed iterative solutions for the variables ${{{{u}}}}({{{{x}}}}, {{{{t}}}})$ and the non-linear term can be expressed as follows:
\begin{equation}\nonumber
\begin{split}
  {{{{{u}}}}}({{{{x}}}},{{{{t}}}})&:=\sum_{j=0}^{\infty}{{{{{u}}}}}_j({{{{x}}}},{{{{t}}}}),\\
  N({{{{{u}}}})({{{{x}}}},{{{{t}}}})}&:={{{{{u}}}}({{{{x}}}},{{{{t}}}})}{\frac {\partial}{{\partial}{{{{x}}}}}}{{{{{u}}}}({{{{x}}}},{{{{t}}}})}=N({{{{{u}}}}_0})+ \sum_{{j}=0}^\infty \bigg[{N} \bigg(\sum_{i=0}^{{j}} {{{{u}}}}_i({{{{x}}}},{{{{t}}}})\bigg)-{N} \bigg(\sum_{i=0}^{{j-1} }{{{{u}}}}_i({{{{x}}}},{{{{t}}}})\bigg)\bigg],
  \end{split}
\end{equation}
where $N(\cdot)$ is the non-linear part as given by \eqref{s4} and $u_j$ will be defined later. After executing the decomposition procedure, equation \eqref{s51} can be expressed in the form:
\begin{equation}\label{qia1}
\begin{split}
\sum_{j=0}^{\infty}{{{{{u}}}}}_j({{{{x}}}},{{{{t}}}})&={{{{{u}}}}}({{{{x}}}},0)+\mathcal{A}^{-1}\left[\frac{1}{s^{{{{\alpha}}}}}\mathcal{A}\left[ f({{{{x}}}},{{{{t}}}})\right]\right]\\
&\qquad+\mathcal{A}^{-1}\frac{1}{s^{{{{\alpha}}}}}\mathcal{A}\begin{cases}\int_{0}^{{{{{t}}}}}\!
    \left( {{{{{t}}}}}-{q} \right) ^{{{{{\alpha}}}}-1}\,{\frac {\partial^2}{{\partial} {{{{x}}}}^2 }}
\left(\sum_{j=0}^{\infty}{{{{{u}}}}}_j({{{{x}}}},{{{{t}}}})\right){\rm d}q-N({{{{{u}}}}_0})\\+ \sum_{{j}=0}^\infty \bigg[{N} \bigg(\sum_{i=0}^{j } {{{{u}}}}_i({{{{x}}}},{{{{t}}}})\bigg)-{N} \bigg(\sum_{i=0}^{{j}-1}{{{{u}}}}_i({{{{x}}}}.{{{{t}}}})\bigg)\bigg].
\end{cases}
\end{split}
\end{equation}

Based on equation \eqref{qia1}, we are ready to provide the recursive IATM algorithm by defining $u_j$ for $j\ge0$. First of all, for $u_0$ and $u_1$, they can be defined respectively as follows:
\begin{equation}\label{ppq}
\begin{split}
&{{{{{u}}}}}_0({{{{x}}}},{{{{t}}}}):={{{{{u}}}}}({{{{x}}}},0)+\mathcal{A}^{-1}\left[\frac{1}{s^{{{{\alpha}}}}}\mathcal{A}\left[f({{{{x}}}},{{{{t}}}})\right]\right],
\end{split}
\end{equation}
\begin{equation}\label{x}
  {{{{{u}}}}}_1({{{{x}}}},{{{{t}}}}):=\mathcal{A}^{-1}\left[\frac{1}{s^{{{{\alpha}}}}}\mathcal{A}\left[\int_{0}^{{{{{t}}}}}\!
    \left( {{{{{t}}}}}-{q} \right) ^{{{{{\alpha}}}}-1}\,{\frac {\partial^2}{{\partial} {{{{x}}}}^2 }}
\left( {{{{{u}}}}}_0({{{{x}}}},{{{{t}}}})\right){\rm d}q-N({{{{{u}}}}_0})\right]\right].
\end{equation}
 For the cases when $j\geq2$, $u_{j}$ can be recursively defined by
\begin{equation}\label{p111}
\begin{split}
{{{{{u}}}}}_{j}({{{{x}}}},{{{{t}}}})=\mathcal{A}^{-1}\left[\frac{1}{s^{{{{\alpha}}}}}\mathcal{A}\left[\int_{0}^{{{{{t}}}}}\!
    \left( {{{{{t}}}}}-{q} \right) ^{{{{{\alpha}}}}-1}\,{\frac {\partial^2}{{\partial} {{{{x}}}}^2 }}
{{{{{u}}}}}_{j-1}({{{{x}}}},{{{{t}}}}){\rm d}q- \bigg[{N} \bigg(\sum_{i=0}^{j-1} {{{{u}}}}_i({{{{x}}}},{{{{t}}}})\bigg)-{N} \bigg(\sum_{i=0}^{{j}-2}{{{{u}}}}_i({{{{x}}}},{{{{t}}}})\bigg)\bigg]\right]\right].
\end{split}
\end{equation}
{Hence the approximate solution $u(x,t)$ can be obtained as a series form $u=\sum_{j=0}^\infty u_j$ with $u_j$ being given by equations \eqref{ppq}, \eqref{x} and \eqref{p111}. The  IATM  is highly effective for addressing a variety of fractional order partial differential equations and fractional integro-differential equations. For those interested in a deeper understanding of its applications, we recommend  references \cite{M, qas1, khan2025efficient}.}
\section{Numerical Implementation for IATM}\label{sec4}
{In this section, we test the validity of IATM by applying the method to various FPIDEs. We aim at considering FPIDEs with different initial conditions and source terms. Similar FPIDEs have been previously investigated by a number of researchers, for example Guo et al. \cite{ex1}, Rawani et al. \cite{rawani2023novel}, T. Akram et al.  \cite{ex12} and the references cited therein. The most up-to-date one was by Khan et al. \cite{LADM} in which the authors achieved higher accuracy by utilizing the Laplace Adomian Decomposition method (LADM).

Throughout this section, all the FPIDEs are posted on $(x,t)\in[0,5]\times[0,1]$ with order $\alpha\in(0,1]$ on the fractional time derivative $D^\alpha_t$. For each problem, the approximate solution is provided by IATM with the iterations being given by equations \eqref{ppq}, \eqref{x} and \eqref{p111}.

We focus on the following non-linear FPIDE
\begin{equation}\label{ss3}
\begin{split}
{{D}^{{{{\alpha}}}}_{{{{{t}}}}}{{{{{u}}}}({{{{x}}}},{{{{t}}}})}}+{{{{{u}}}}({{{{x}}}},{{{{t}}}})}{\frac {\partial}{{\partial}{{{{x}}}}}}{{{{{u}}}}({{{{x}}}},{{{{t}}}})}=\int_{0}^{{{{{t}}}}}\!
   \left( {{{{{t}}}}}-{q} \right) ^{{{{{\alpha}}}}-1}\,{\frac {\partial^2}{{\partial} {{{{x}}}}^2 }}
{{{{{u}}}}({{{{x}}}},{{{{t}}}})} {\rm d}q+f({{{{x}}}},{{{{t}}}})
\end{split}
\end{equation}
with different initial conditions $u(x,0)$ and source terms $f(x,t)$ which will be provided later. Equation \eqref{ss3} can be used for modelling physical phenomena which involve heat flow in materials with memory and phenomena associated with linear viscoelastic mechanics \cite{ex1, rawani2023novel}. The non-linear term $N$ in equation \eqref{ss3} is given by
\begin{equation}\label{non-linear term main}
\begin{split}
  N({{{{{u}}}})({{{{x}}}},{{{{t}}}})}&:={{{{{u}}}}({{{{x}}}},{{{{t}}}})}{\frac {\partial}{{\partial}{{{{x}}}}}}{{{{{u}}}}({{{{x}}}},{{{{t}}}})}.
  \end{split}
\end{equation}
The Adomian polynomials $\mathbf{P}_j $ \cite{LADM} and the Daftardar--Jafari polynomials $\mathbf{J}_j$ are defined respectively as follows:
\begin{equation*}
\begin{split}
&\textnormal {\underline{Adomian Polynomials}} \ \ \ \ \ \  \ \ \ \ \ \  \ \ \  \ \ \      \ \ \ \ \ \   \ \ \   \ \ \   \ \ \   \ \ \   \ \ \ \textnormal {\underline{Daftardar--Jafari polynomials}} \\
\mathbf{P}_{0} &:= {{{{{u}}}}_{0}({{{{x}}}},{{{{t}}}})}{\frac {\partial}{{\partial}{{{{x}}}}}}{{{{{u}}}}_{0}({{{{x}}}},{{{{t}}}})}, \ \ \ \ \ \ \ \ \  \ \ \ \ \ \  \ \ \  \ \ \      \ \ \   \ \ \   \ \ \   \ \ \   \ \ \   \ \ \  \mathbf{J}_{0} := {{{{{u}}}}_{0}({{{{x}}}},{{{{t}}}})}{\frac {\partial}{{\partial}{{{{x}}}}}}{{{{{u}}}}_{0}({{{{x}}}},{{{{t}}}})},\\
\mathbf{P}_{1} &:= u_{{1}} \left( x,t \right) {\frac {\partial }{\partial x}}u_{{0}}
 \left( x,t \right) +u_{{0}} \left( x,t \right) {\frac {\partial }{
\partial x}}u_{{1}} \left( x,t \right)
, \ \ \  \ \ \    \ \ \     \mathbf{J}_{1} := \begin{cases}u_{{0}} \left( x,t \right) {\frac {\partial }{\partial x}}u_{{1}}
 \left( x,t \right) +u_{{1}} \left( x,t \right) {\frac {\partial }{
\partial x}}u_{{0}} \left( x,t \right) \\    +u_{{1}} \left( x,t \right) {
\frac {\partial }{\partial x}}u_{{1}} \left( x,t \right)
, \end{cases}\\
\mathbf{P}_{2} &:=\begin{cases} u_{{2}} \left( x,t \right) {\frac {\partial }{\partial x}}u_{{0}}
 \left( x,t \right) +u_{{1}} \left( x,t \right) {\frac {\partial }{
\partial x}}u_{{1}} \left( x,t \right) \\+u_{{0}} \left( x,t \right) {
\frac {\partial }{\partial x}}u_{{2}} \left( x,t \right)
,\end{cases} \ \ \   \ \ \   \              \mathbf{J}_{2} :=\begin{cases}u_{{0}} \left( x,t \right) {\frac {\partial }{\partial x}}u_{{2}}
 \left( x,t \right) +u_{{1}} \left( x,t \right) {\frac {\partial }{
\partial x}}u_{{2}} \left( x,t \right)\\ +u_{{2}} \left( x,t \right) {
\frac {\partial }{\partial x}}u_{{0}} \left( x,t \right) \\+u_{{2}}
 \left( x,t \right) {\frac {\partial }{\partial x}}u_{{1}} \left( x,t
 \right) +u_{{2}} \left( x,t \right) {\frac {\partial }{\partial x}}u_
{{2}} \left( x,t \right)
.     \end{cases} \\
 \vdots
\end{split}
\end{equation*}

Equation \eqref{ss3} will be equipped with different initial conditions $u(x,0)$ and source terms $f(x,t)$, which are listed in Problem~\ref{TE2}, Problem~\ref{TE4} and Problem~\ref{TE1} below.}

\begin{prb}\label{TE2}
We consider the FPIDE \eqref{ss3} with source term $f(x,t)$
\begin{equation}\label{st1}
\begin{split}
f({{{{x}}}},{{{{t}}}})&=\pi\, \left( {\frac {\pi\,{{{{{t}}}}}^{{{{{\alpha}}}}}}{{{{{\alpha}}}}}}-4\,{{{{{t}}}}}^{3}\cos
 \left( 2\,\pi\,{{{{x}}}} \right)  \right) \sin \left( \pi\,{{{{x}}}} \right) +\pi/2-
12\,{\frac {{{{{{t}}}}}^{3-{{{{\alpha}}}}}}{\Gamma \left( 4-{{{{\alpha}}}} \right) }}-2\,\pi\,
{{{{{t}}}}}^{3}\cos \left( \pi\,{{{{x}}}} \right)\\
& \qquad -48\,{\frac {{\pi}^{2}\Gamma \left(
{{{{\alpha}}}} \right) {{{{{t}}}}}^{3+{{{{\alpha}}}}}}{\Gamma \left( 4+{{{{\alpha}}}} \right) }}+8\,
\pi\,{{{{{t}}}}}^{6}\cos \left( 2\,\pi\,{{{{x}}}} \right) \sin \left( 2\,\pi\,{{{{x}}}}
 \right)
 \end{split}
\end{equation}
and initial condition $u(x,0)$
\begin{equation}\label{ic1}
{{{{u}}}}({{{{x}}}},0)= \sin(\pi {{{{{x}}}}}).
\end{equation}
The exact solution to equations \eqref{ss3}, \eqref{st1} and \eqref{ic1} is
\begin{equation*}
  {{{{u}}}}({{{{x}}}},{{{{t}}}})=\sin \left( \pi\,{{{{x}}}} \right) -2\,{{{{{t}}}}}^{3}\sin \left( 2\,\pi\,{{{{x}}}} \right).
\end{equation*}
{Table \ref{table1} shows the numerical simulations for the approximate solution and the exact solution respectively for Problem~\ref{TE2}. Figure \ref{figure1} shows the comparison solution plots for Problem~\ref{TE2} for various values of $\alpha$, and Figure \ref{figure2} shows the comparison plots between approximate solution and exact solution at $\alpha=1$.}
%Table 1
 \begin{table}[H]
				\centering
				\caption{{Numerical simulations across various iterations, time levels, and spatial domains for the Problem  \ref{TE2} at $\alpha=1$.}}\label{table1}
				\resizebox{\textwidth}{!}{
					\begin{tabular}{ccccccc}
						\toprule
					 	{${{{t}}}$}& ${{x}}$& {\thead{ Approximate \\    solution}} &{\thead{ Exact\\ solution}}  & {\thead{ Absolute error\\ of LADM \cite{LADM}\\at 2nd iteration}}&  {\thead{ Absolute error\\  of our approach \\ at 2nd iteration}} &{\thead  { Absolute error\\  of our approach \\ at 3rd iteration} }\\ %[0.1ex]
						\midrule
\num{0.001}&\num{0.2}&\num{0.587785250}&\num{0.587785252}&\num{2.82E-06}&\num{1.55E-09}&\num{1.27E-09}\\
 \num{0.001}&\num{0.4}&\num{0.951056515}&\num{0.951056518}&\num{3.33E-06}&\num{1.33E-08}&\num{2.80E-09}\\
 \num{0.001}&\num{0.6}&\num{0.951056517}&\num{0.951056518}&\num{3.36E-06}&\num{9.89E-09}&\num{4.54E-10}\\
 \num{0.001}&\num{0.8}&\num{0.587785254}&\num{0.587785252}&\num{2.77E-06}&\num{2.97E-09}&\num{2.53E-09}\\
 \num{0.001}&\num{1.0}&\num{-0.000000000}&\num{-0.000000002}&\num{2.00E-09}&\num{2.00E-09}&\num{2.00E-09}\\
 \num{0.003}&\num{0.2}&\num{0.587785202}&\num{0.587785236}&\num{2.59E-05}&\num{3.65E-08}&\num{3.36E-08}\\
 \num{0.003}&\num{0.4}&\num{0.951056484}&\num{0.951056560}&\num{2.97E-05}&\num{3.63E-07}&\num{7.62E-08}\\
 \num{0.003}&\num{0.6}&\num{0.951056547}&\num{0.951056560}&\num{3.05E-05}&\num{2.63E-07}&\num{1.29E-08}\\
 \num{0.003}&\num{0.8}&\num{0.587785305}&\num{0.587785236}&\num{2.44E-05}&\num{8.55E-08}&\num{6.91E-08}\\
 \num{0.003}&\num{1.0}&\num{-0.000000000}&\num{-0.000000054}&\num{5.40E-08}&\num{5.40E-08}&\num{5.40E-08}\\
 \num{0.005}&\num{0.2}&\num{0.587785022}&\num{0.587785175}&\num{7.33E-05}&\num{1.44E-07}&\num{1.53E-07}\\
 \num{0.005}&\num{0.4}&\num{0.951056364}&\num{0.951056719}&\num{8.18E-05}&\num{1.70E-06}&\num{3.55E-07}\\
 \num{0.005}&\num{0.6}&\num{0.951056655}&\num{0.951056719}&\num{8.53E-05}&\num{1.20E-06}&\num{6.34E-08}\\
 \num{0.005}&\num{0.8}&\num{0.587785499}&\num{0.587785175}&\num{6.63E-05}&\num{4.19E-07}&\num{3.24E-07}\\
 \num{0.005}&\num{1.0}&\num{-0.000000000}&\num{-0.000000250}&\num{2.50E-07}&\num{2.50E-07}&\num{2.50E-07}\\
 \num{0.007}&\num{0.2}&\num{0.587784630}&\num{0.587785040}&\num{1.46E-04}&\num{3.27E-07}&\num{4.10E-07}\\
 \num{0.007}&\num{0.4}&\num{0.951056094}&\num{0.951057071}&\num{1.59E-04}&\num{4.71E-06}&\num{9.77E-07}\\
 \num{0.007}&\num{0.6}&\num{0.951056887}&\num{0.951057071}&\num{1.68E-04}&\num{3.23E-06}&\num{1.85E-07}\\
 \num{0.007}&\num{0.8}&\num{0.587785938}&\num{0.587785040}&\num{1.27E-04}&\num{1.21E-06}&\num{8.98E-07}\\
 \num{0.007}&\num{1.0}&\num{-0.000000000}&\num{-0.000000686}&\num{6.86E-07}&\num{6.86E-07}&\num{6.86E-07}\\
 \num{0.009}&\num{0.2}&\num{0.587783946}&\num{0.587784802}&\num{2.46E-04}&\num{5.49E-07}&\num{8.56E-07}\\
 \num{0.009}&\num{0.4}&\num{0.951055613}&\num{0.951057696}&\num{2.60E-04}&\num{1.01E-05}&\num{2.08E-06}\\
 \num{0.009}&\num{0.6}&\num{0.951057278}&\num{0.951057696}&\num{2.80E-04}&\num{6.76E-06}&\num{4.18E-07}\\
 \num{0.009}&\num{0.8}&\num{0.587786731}&\num{0.587784802}&\num{2.06E-04}&\num{2.72E-06}&\num{1.93E-06}\\
 \num{0.009}&\num{1.0}&\num{-0.000000000}&\num{-0.000001458}&\num{1.46E-06}&\num{1.46E-06}&\num{1.46E-06 }\\
 \bottomrule
				\end{tabular}}
			\end{table}
%Figure 1%			
\begin{figure}[H]
	\centering
	\begin{subfigure}{0.31\textwidth}
		\includegraphics[width=6.2cm]{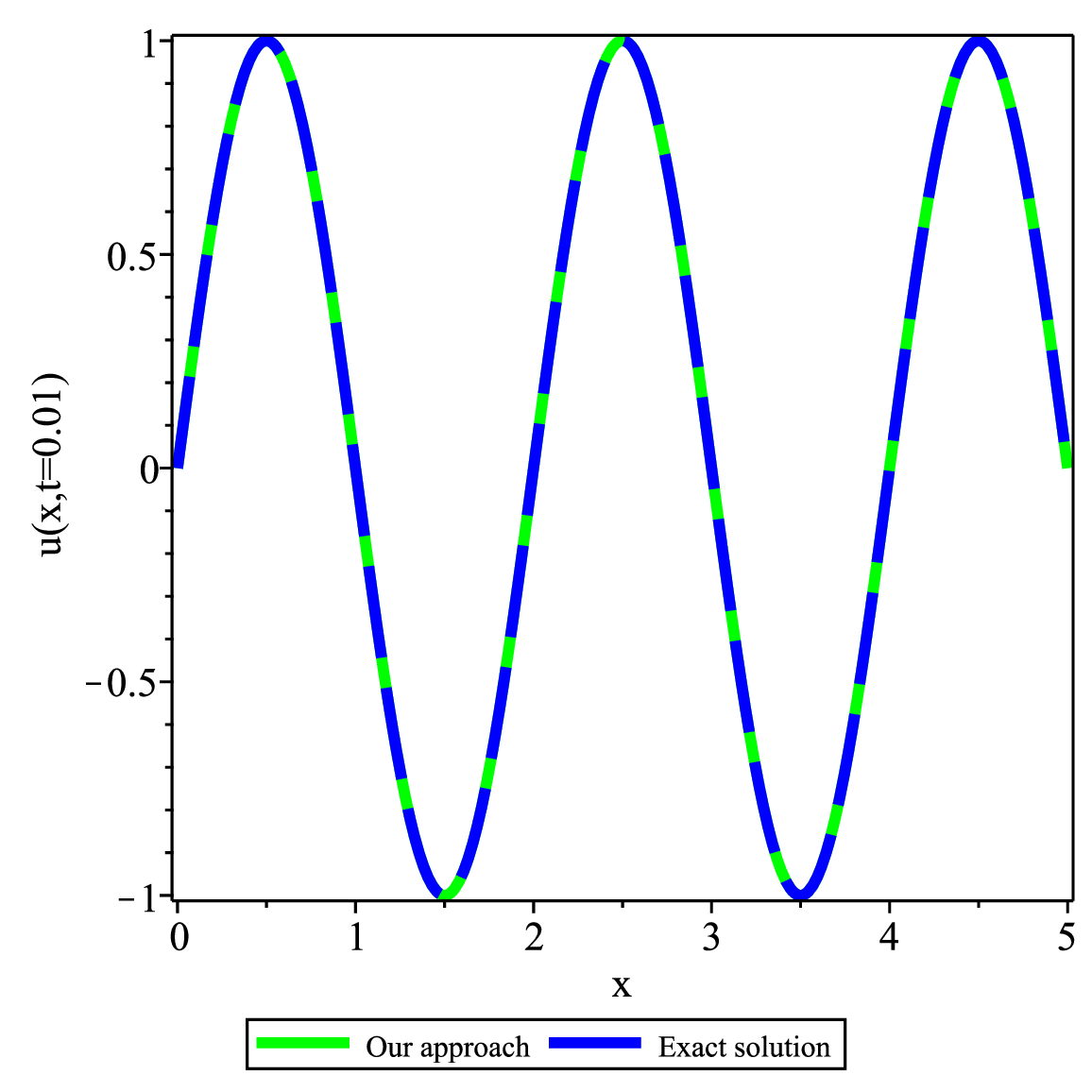}
		\caption{Exact and approximate solution  at $\alpha=1$}
		\label{fig:noise_stabilize_norm}
	\end{subfigure}\hspace{6.6em}%\hfill
	\begin{subfigure}{0.31\textwidth}
		\includegraphics[width=6.2cm]{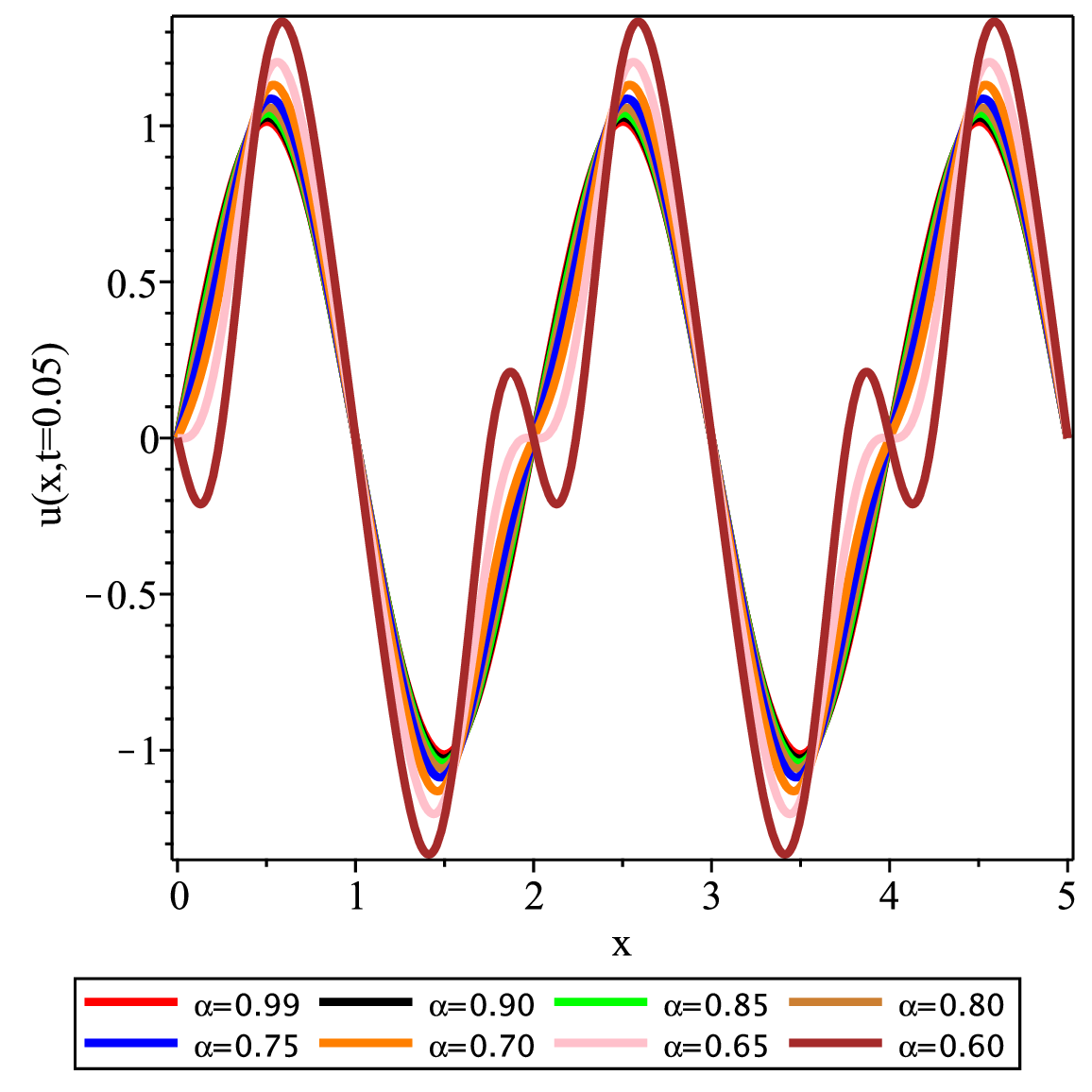}
		\caption{Solution plot for different   $\alpha$}
		\label{figure1x}
	\end{subfigure}\hspace{5.1em}
	\caption{Comparison solution plots  { up to three terms approximation for Problem \ref{TE2}.}\label{figure1}}
\end{figure}
%Figure 2 at alpha=1
\begin{figure}[H]
	\centering 
	\begin{subfigure}{0.31\textwidth}
		\includegraphics[width=6.2cm]{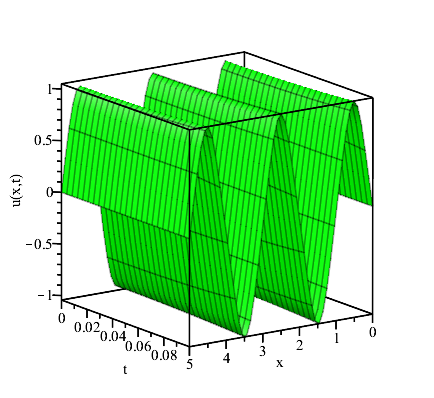}
		\caption{3D approximate  solution}
		\label{fig:noise_stabilize_u}
	\end{subfigure}\hspace{6.9em}
	%\hfill
	\begin{subfigure}{0.31\textwidth}
		\includegraphics[width=6.2cm]{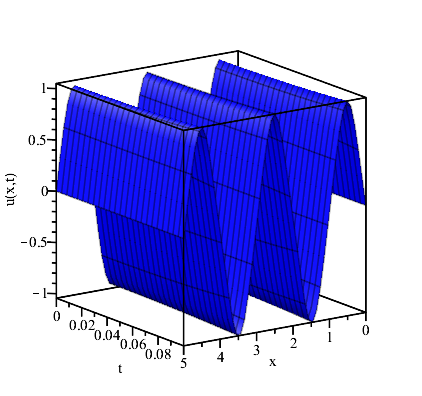}
		\caption{3D exact solution }
		\label{figa1}
	\end{subfigure}\hspace{0.1em}
	\caption{Comparison solution plots { up to three terms approximation for Problem \ref{TE2} at $\alpha=1$.}}\label{figure2}
\end{figure}
\end{prb}

\begin{prb}\label{TE4}
Next, we consider the FPIDE \eqref{ss3} with source term $f(x,t)$
\begin{equation}\label{st2}
\begin{split}
f({{{{x}}}},{{{{t}}}})&= \left( 6\,{\frac {{t}^{3-\alpha}}{\Gamma \left( 4-\alpha \right) }}+
\pi\,{t}^{6}\cos \left( \pi\,x \right) +6\,{\frac {{\pi}^{2}\Gamma
 \left( \alpha \right) {t}^{3+\alpha}}{\Gamma \left( 4+\alpha \right)
}} \right) \sin \left( \pi\,x \right)
 \end{split}
\end{equation}
and initial condition $u(x,0)$
\begin{equation}\label{ic2}
{{{{u}}}}({{{{x}}}},0)= 0.
\end{equation}
The exact solution to equations \eqref{ss3}, \eqref{st2} and \eqref{ic2} is
\begin{equation*}
  {{{{u}}}}({{{{x}}}},{{{{t}}}})={t}^{3}\sin \left( \pi\,x \right).
\end{equation*}
{Table \ref{table2} shows the numerical simulations for the approximate solution and the exact solution respectively for Problem~\ref{TE4}. Figure \ref{figure3} shows the comparison solution plots for Problem~\ref{TE4} for various values of $\alpha$, and Figure \ref{figure4} shows the comparison plots between approximate solution and exact solution at $\alpha=1$.}
%Figure 3			
\begin{figure}[H]
	\centering
	\begin{subfigure}{0.31\textwidth}
		\includegraphics[width=6.2cm]{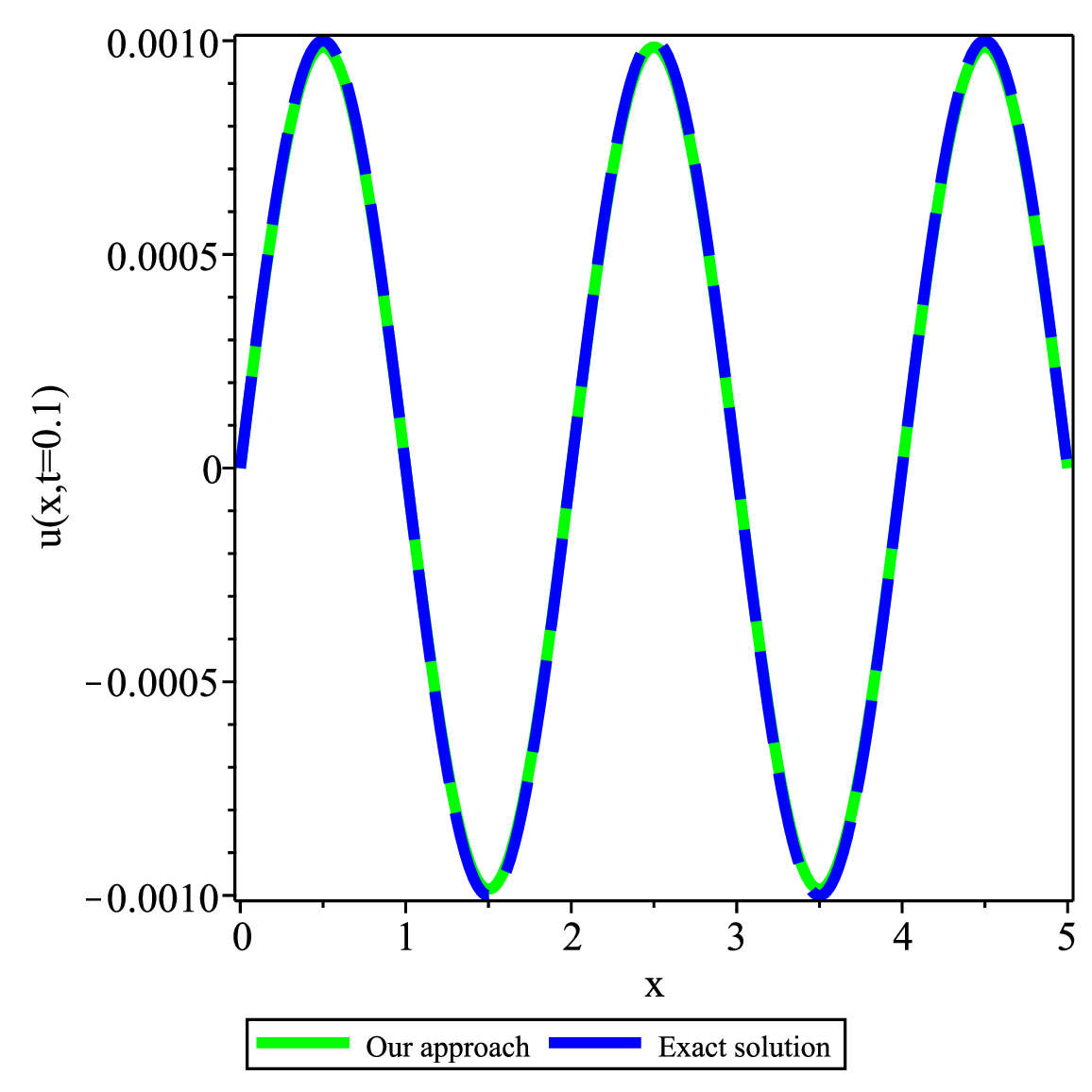}
		\caption{Exact and approximate solution  at $\alpha=1$}
		\label{fig:noise_stabilize_norm}
	\end{subfigure}\hspace{6.6em}%\hfill
	\begin{subfigure}{0.31\textwidth}
		\includegraphics[width=6.2cm]{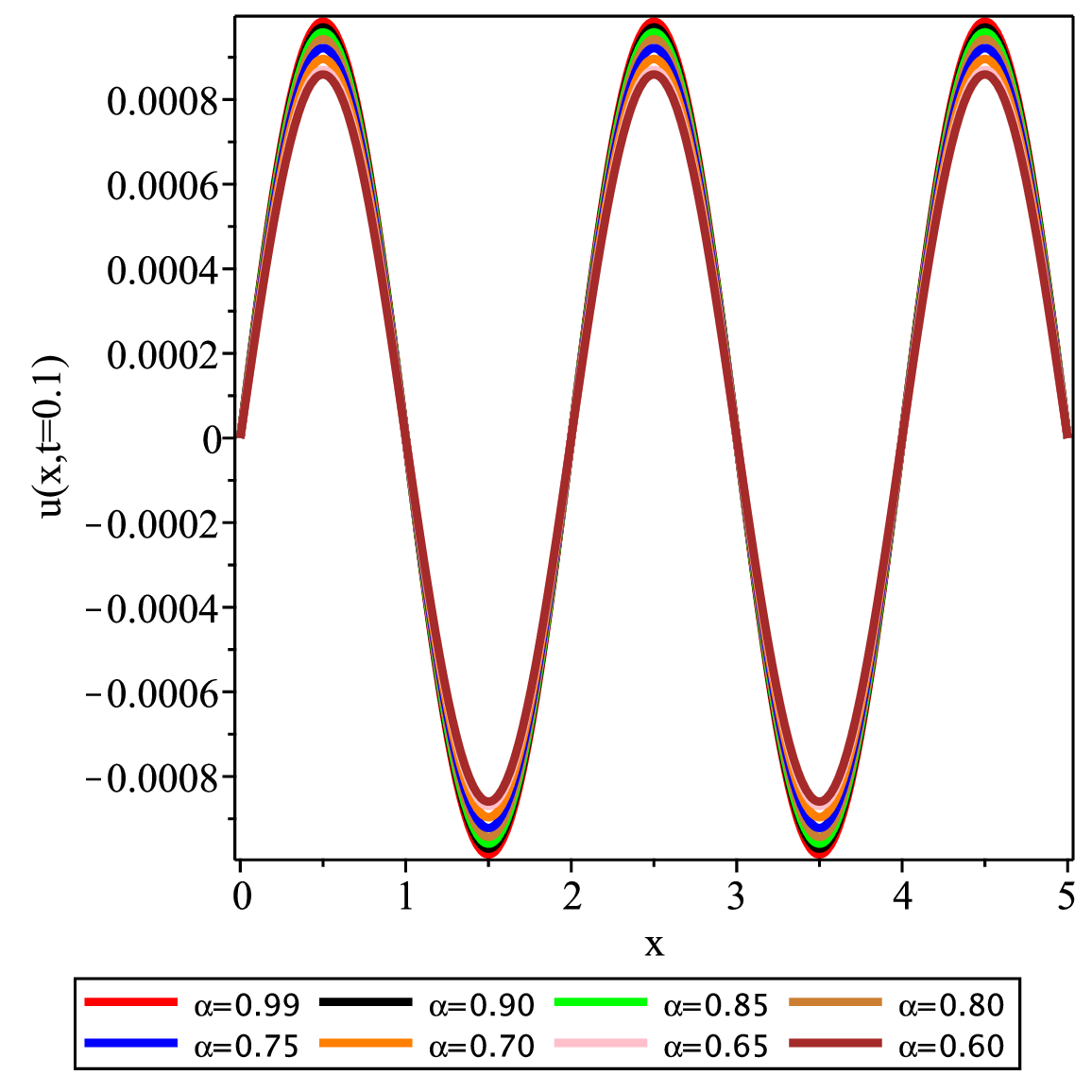}
		\caption{Solution plot for different   $\alpha$}
		\label{figure2x}
	\end{subfigure}\hspace{6.1em}
	\caption{Comparison solution plots  { up to three terms approximation for Problem \ref{TE4}.}}\label{figure3}
\end{figure}
%Figure 4 at alpha=1
			\begin{figure}[H]
	\centering
	\begin{subfigure}{0.31\textwidth}
		\includegraphics[width=6.2cm]{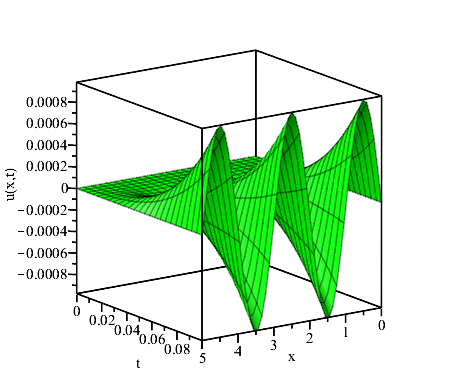}
		\caption{3D approximate  solution}
		\label{ewd}
	\end{subfigure}\hspace{6.9em}
	%\hfill
	\begin{subfigure}{0.31\textwidth}
		\includegraphics[width=6.2cm]{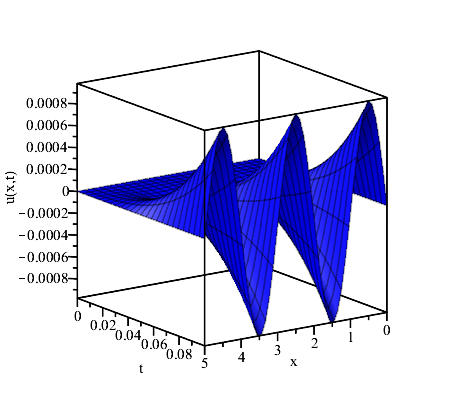}
		\caption{3D exact solution }
		\label{figa2}
	\end{subfigure}\hspace{0.9em}
	\caption{Comparison solution plots { up to three terms approximation for Problem \ref{TE4} at $\alpha=1$.}}\label{figure4}
\end{figure}

%Table 2
 \begin{table}[H]
				\centering
				\caption{{Numerical simulations across various iterations, time levels, and spatial domains for the Problem  \ref{TE4} at $\alpha=1$.}}\label{table2}
				\resizebox{\textwidth}{!}{
					\begin{tabular}{ccccccc}
						\toprule
					 	{${{{t}}}$}& ${{x}}$& {\thead{ Approximate \\    solution}} &{\thead{ Exact\\ solution}}  & {\thead{ Absolute error\\ at 1st iteration}}&  {\thead{ Absolute error\\ at 2nd iteration}} &{\thead  { Absolute error\\ at 3rd iteration} }\\ %[0.1ex]
						\midrule
\num{0.01} & \num{0.1} & \num{0.000000309} & \num{0.000000309} & \num{6.10E-11} & \num{8.60E-15} & \num{9.44E-19} \\
\num{0.01} & \num{0.3} & \num{0.000000809} & \num{0.000000809} & \num{1.60E-10} & \num{2.25E-14} & \num{2.47E-18} \\
\num{0.01} & \num{0.5} & \num{0.000001000} & \num{0.000001000} & \num{1.97E-10} & \num{2.78E-14} & \num{3.05E-18} \\
\num{0.01} & \num{0.7} & \num{0.000000809} & \num{0.000000809} & \num{1.60E-10} & \num{2.25E-14} & \num{2.47E-18} \\
\num{0.01} & \num{0.9} & \num{0.000000309} & \num{0.000000309} & \num{6.10E-11} & \num{8.60E-15} & \num{9.43E-19} \\
\num{0.03} & \num{0.1} & \num{0.000008343} & \num{0.000008343} & \num{1.48E-08} & \num{1.88E-11} & \num{1.86E-14} \\
\num{0.03} & \num{0.3} & \num{0.000021843} & \num{0.000021843} & \num{3.88E-08} & \num{4.93E-11} & \num{4.87E-14} \\
\num{0.03} & \num{0.5} & \num{0.000027000} & \num{0.000027000} & \num{4.80E-08} & \num{6.09E-11} & \num{6.01E-14} \\
\num{0.03} & \num{0.7} & \num{0.000021843} & \num{0.000021843} & \num{3.88E-08} & \num{4.92E-11} & \num{4.85E-14} \\
\num{0.03} & \num{0.9} & \num{0.000008343} & \num{0.000008343} & \num{1.48E-08} & \num{1.88E-11} & \num{1.85E-14} \\
\num{0.05} & \num{0.1} & \num{0.000038627} & \num{0.000038627} & \num{1.91E-07} & \num{6.74E-10} & \num{1.86E-12} \\
\num{0.05} & \num{0.3} & \num{0.000101127} & \num{0.000101127} & \num{4.99E-07} & \num{1.76E-09} & \num{4.85E-12} \\
\num{0.05} & \num{0.5} & \num{0.000125000} & \num{0.000125000} & \num{6.17E-07} & \num{2.17E-09} & \num{5.96E-12} \\
\num{0.05} & \num{0.7} & \num{0.000101127} & \num{0.000101127} & \num{4.99E-07} & \num{1.76E-09} & \num{4.79E-12} \\
\num{0.05} & \num{0.9} & \num{0.000038627} & \num{0.000038627} & \num{1.91E-07} & \num{6.70E-10} & \num{1.82E-12} \\
\num{0.07} & \num{0.1} & \num{0.000105993} & \num{0.000105993} & \num{1.03E-06} & \num{7.12E-09} & \num{3.89E-11} \\
\num{0.07} & \num{0.3} & \num{0.000277493} & \num{0.000277493} & \num{2.69E-06} & \num{1.86E-08} & \num{1.01E-10} \\
\num{0.07} & \num{0.5} & \num{0.000343000} & \num{0.000343000} & \num{3.32E-06} & \num{2.29E-08} & \num{1.23E-10} \\
\num{0.07} & \num{0.7} & \num{0.000277493} & \num{0.000277493} & \num{2.68E-06} & \num{1.85E-08} & \num{9.84E-11} \\
\num{0.07} & \num{0.9} & \num{0.000105993} & \num{0.000105993} & \num{1.02E-06} & \num{7.04E-09} & \num{3.73E-11} \\
\num{0.09} & \num{0.1} & \num{0.000225274} & \num{0.000225273} & \num{3.61E-06} & \num{4.15E-08} & \num{3.78E-10} \\
\num{0.09} & \num{0.3} & \num{0.000589774} & \num{0.000589773} & \num{9.44E-06} & \num{1.08E-07} & \num{9.77E-10} \\
\num{0.09} & \num{0.5} & \num{0.000729001} & \num{0.000729000} & \num{1.17E-05} & \num{1.33E-07} & \num{1.18E-09} \\
\num{0.09} & \num{0.7} & \num{0.000589774} & \num{0.000589773} & \num{9.42E-06} & \num{1.07E-07} & \num{9.36E-10} \\
\num{0.09} & \num{0.9} & \num{0.000225274} & \num{0.000225273} & \num{3.60E-06} & \num{4.08E-08} & \num{3.53E-10} \\
 \bottomrule
				\end{tabular}}
			\end{table}

\end{prb}

\begin{prb}\label{TE1}
Finally, we consider the FPIDE \eqref{ss3} with source term $f(x,t)$
\begin{equation}\label{st3}
\begin{split}
f({{{{x}}}},{{{{t}}}})&={\frac {\Gamma \left( 7/2 \right) {{{{{t}}}}}^{5/2-{{{{\alpha}}}}}{{{{{x}}}}}^{2} \left( 1-{{{{{x}}}}}
 \right) ^{2}}{\Gamma \left( 7/2-{{{{\alpha}}}} \right) }}-2\, \left( {\frac {
{{{{{t}}}}}^{{{{{\alpha}}}}}}{{{{{\alpha}}}}}}+\,{\frac {\Gamma \left( 7/2 \right) \Gamma
 \left( {{{{\alpha}}}} \right) {{{{{t}}}}}^{5/2+{{{{\alpha}}}}}}{\Gamma \left( 7/2+{{{{\alpha}}}}
 \right) }} \right)  \left( 6\,{{{{{x}}}}}^{2}-6\,{{{{{x}}}}}+1 \right) \\&+2\, \left( 1+{{{{{t}}}}
}^{5/2} \right) ^{2} \left( 1-2\,{{{{{x}}}}}\right) {{{{{x}}}}}^{3} \left( 1-{{{{{x}}}}}\right)
^{3}
 \end{split}
\end{equation}
and initial condition $u(x,0)$
\begin{equation}\label{ic3}
{{{{u}}}}({{{{x}}}},0)= {{{{{x}}}}}^{2} \left( 1-{{{{{x}}}}}\right) ^{2}.
\end{equation}
The exact solution solution to equations \eqref{ss3}, \eqref{st3} and \eqref{ic3} is
\begin{equation*}
  {{{{u}}}}({{{{x}}}},{{{{t}}}})=\left( 1+{{{{{t}}}}}^{5/2} \right) {{{{{x}}}}}^{2} \left( 1-{{{{{x}}}}}\right) ^{2}.
\end{equation*}
{Table \ref{table3} shows the numerical simulations for the approximate solution and the exact solution respectively for Problem~\ref{TE1}. Figure \ref{figure5} shows the comparison solution plots for Problem~\ref{TE1} for various values of $\alpha$, and Figure \ref{figure6} shows the comparison plots between approximate solution and exact solution at $\alpha=1$.}
%Table 3
 \begin{table}[H]
				\centering
				\caption{{Numerical simulations across various iterations, time levels, and spatial domains for the Problem  \ref{TE1} at $\alpha=1$.}}\label{table3}
				\resizebox{\textwidth}{!}{
					\begin{tabular}{ccccccc}
						\toprule
					 	{${{{t}}}$}& ${{x}}$& {\thead{ Approximate \\    solution}} &{\thead{ Exact\\ solution}}  & {\thead{ Absolute error\\ at 1st iteration}}&  {\thead{ Absolute error\\ at 2nd iteration}} &{\thead  { Absolute error\\ at 3rd iteration} }\\ %[0.1ex]
						\midrule
 \num{0.001}&\num{0.2}&\num{0.025600001}&\num{0.025600001}&\num{1.05E-09}&\num{1.59E-13}&\num{1.95E-17}\\
 \num{0.001}&\num{0.4}&\num{0.057600002}&\num{0.057600002}&\num{7.58E-10}&\num{3.98E-13}&\num{2.93E-17}\\
 \num{0.001}&\num{0.6}&\num{0.057600002}&\num{0.057600002}&\num{1.09E-09}&\num{3.35E-13}&\num{2.50E-17}\\
 \num{0.001}&\num{0.8}&\num{0.025600001}&\num{0.025600001}&\num{8.99E-10}&\num{6.89E-14}&\num{8.29E-18}\\
 \num{0.001}&\num{1.0}&\num{-0.000000000}&\num{0.000000000}&\num{3.00E-12}&\num{1.92E-14}&\num{6.05E-18}\\
 \num{0.003}&\num{0.2}&\num{0.025600013}&\num{0.025600013}&\num{1.10E-08}&\num{7.10E-12}&\num{4.95E-15}\\
 \num{0.003}&\num{0.4}&\num{0.057600028}&\num{0.057600028}&\num{3.59E-09}&\num{1.19E-11}&\num{7.24E-15}\\
 \num{0.003}&\num{0.6}&\num{0.057600028}&\num{0.057600028}&\num{1.26E-08}&\num{6.77E-12}&\num{5.88E-15}\\
 \num{0.003}&\num{0.8}&\num{0.025600013}&\num{0.025600013}&\num{6.97E-09}&\num{2.22E-13}&\num{2.07E-15}\\
 \num{0.003}&\num{1.0}&\num{-0.000000000}&\num{0.000000000}&\num{2.43E-10}&\num{4.67E-12}&\num{4.24E-15}\\
 \num{0.005}&\num{0.2}&\num{0.025600045}&\num{0.025600045}&\num{3.54E-08}&\num{4.75E-11}&\num{7.51E-14}\\
 \num{0.005}&\num{0.4}&\num{0.057600102}&\num{0.057600102}&\num{4.00E-10}&\num{5.72E-11}&\num{9.63E-14}\\
 \num{0.005}&\num{0.6}&\num{0.057600102}&\num{0.057600102}&\num{4.23E-08}&\num{1.80E-11}&\num{7.01E-14}\\
 \num{0.005}&\num{0.8}&\num{0.025600045}&\num{0.025600045}&\num{1.68E-08}&\num{9.02E-12}&\num{2.24E-14}\\
 \num{0.005}&\num{1.0}&\num{-0.000000000}&\num{0.000000000}&\num{1.88E-09}&\num{6.00E-11}&\num{8.71E-14}\\
 \num{0.007}&\num{0.2}&\num{0.025600105}&\num{0.025600105}&\num{8.02E-08}&\num{1.75E-10}&\num{4.79E-13}\\
 \num{0.007}&\num{0.4}&\num{0.057600236}&\num{0.057600236}&\num{1.92E-08}&\num{1.56E-10}&\num{5.34E-13}\\
 \num{0.007}&\num{0.6}&\num{0.057600236}&\num{0.057600236}&\num{9.58E-08}&\num{5.31E-12}&\num{3.43E-13}\\
 \num{0.007}&\num{0.8}&\num{0.025600105}&\num{0.025600105}&\num{2.93E-08}&\num{4.21E-11}&\num{1.02E-13}\\
 \num{0.007}&\num{1.0}&\num{-0.000000000}&\num{0.000000000}&\num{7.22E-09}&\num{3.23E-10}&\num{6.28E-13}\\
 \num{0.009}&\num{0.2}&\num{0.025600197}&\num{0.025600197}&\num{1.52E-07}&\num{4.76E-10}&\num{1.99E-12}\\
 \num{0.009}&\num{0.4}&\num{0.057600443}&\num{0.057600443}&\num{6.66E-08}&\num{3.13E-10}&\num{1.92E-12}\\
 \num{0.009}&\num{0.6}&\num{0.057600443}&\num{0.057600443}&\num{1.78E-07}&\num{9.84E-11}&\num{1.07E-12}\\
 \num{0.009}&\num{0.8}&\num{0.025600197}&\num{0.025600197}&\num{4.41E-08}&\num{1.17E-10}&\num{3.19E-13}\\
 \num{0.009}&\num{1.0}&\num{-0.000000000}&\num{0.000000000}&\num{1.98E-08}&\num{1.13E-09}&\num{2.71E-12}\\
 \bottomrule
				\end{tabular}}
			\end{table}
%Figure 5
\begin{figure}[H]
	\centering
	\begin{subfigure}{0.31\textwidth}
		\includegraphics[width=6.2cm]{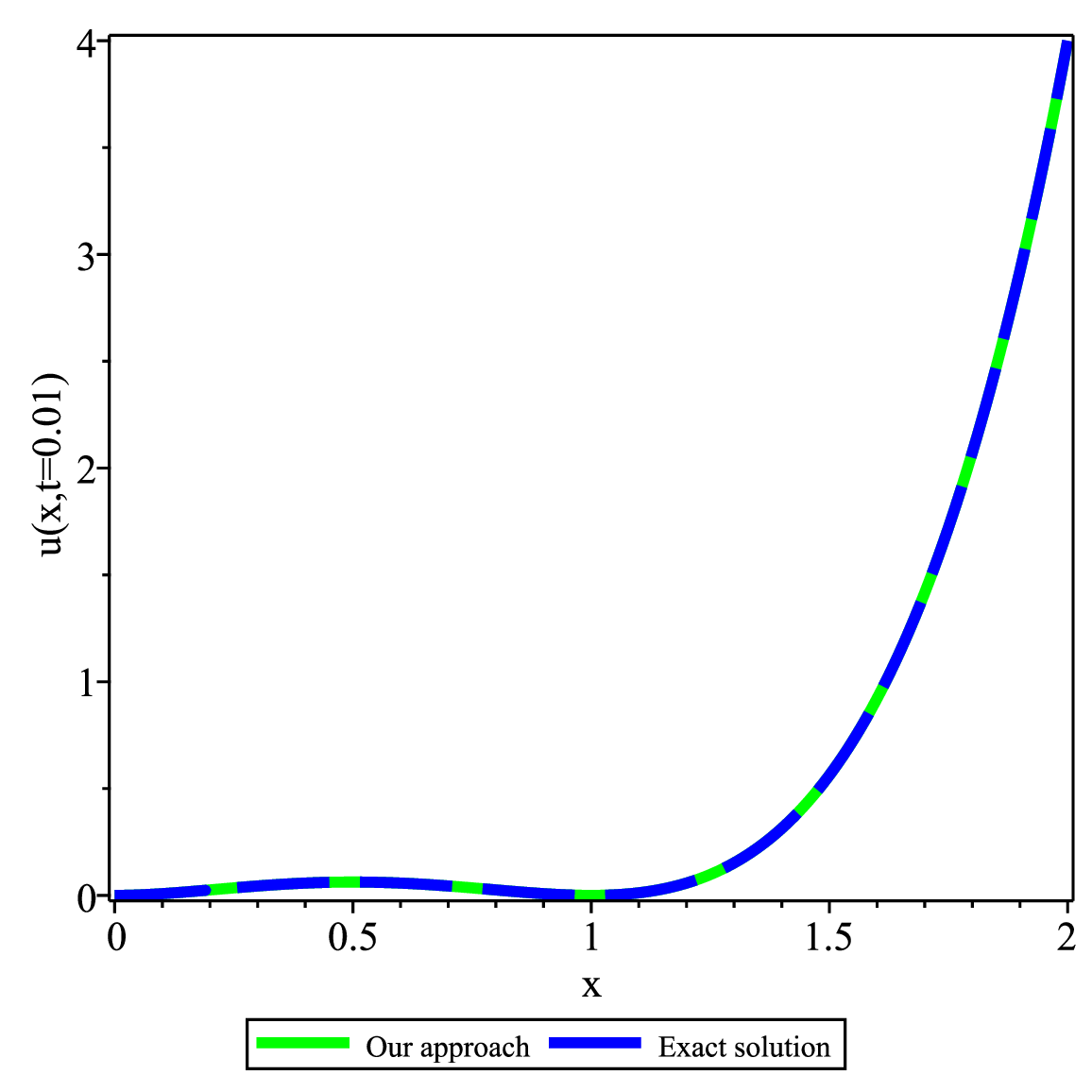}
		\caption{Exact and approximate solution}
		\label{fig:noise_stabilize_norm}
	\end{subfigure}\hspace{8.9em}%\hfill
\begin{subfigure}{0.31\textwidth}
		\includegraphics[width=6.2cm]{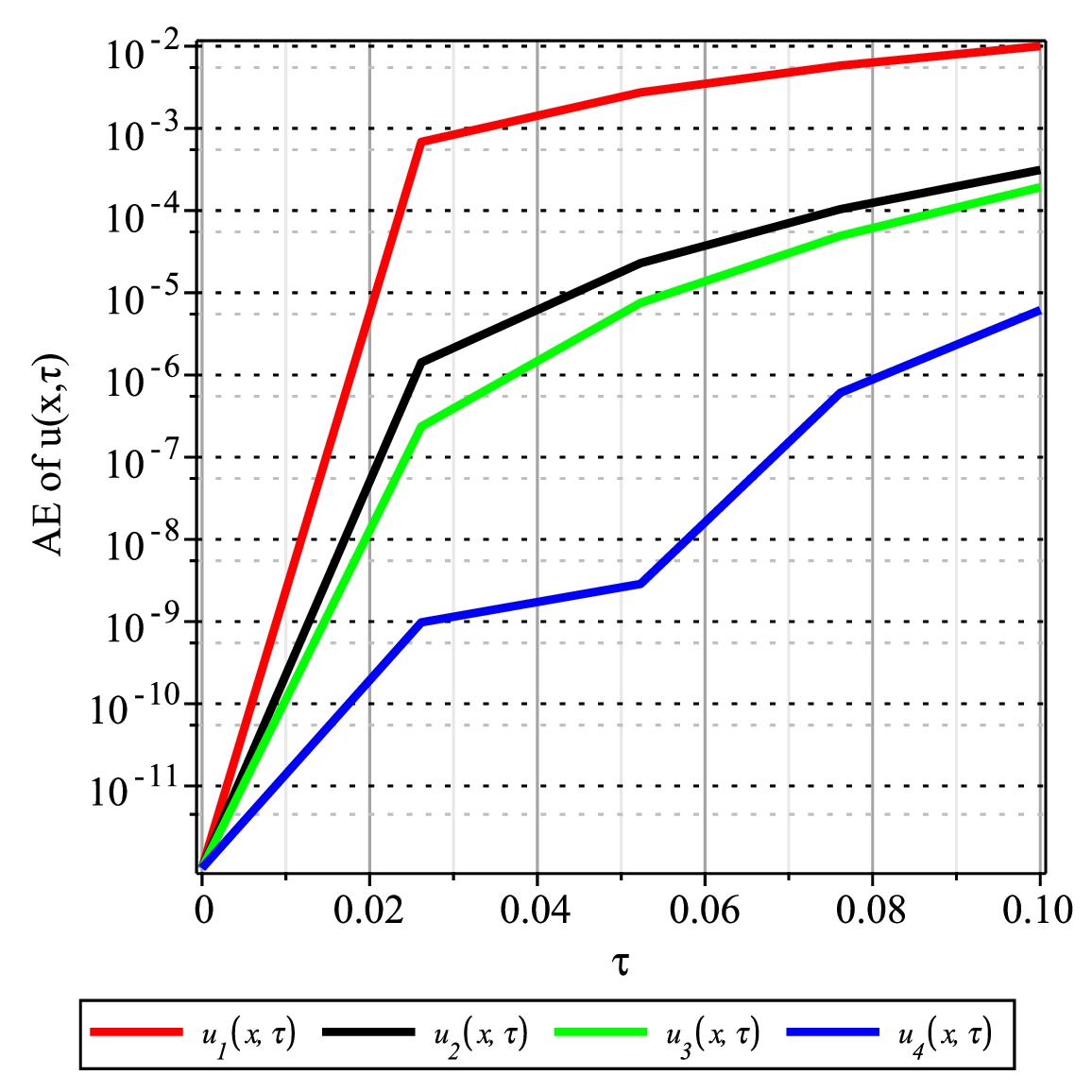}
		\caption{Error for different iterations}
		\label{figure3x}
	\end{subfigure}\hspace{7.1em}
	\caption{Comparison solution plots  { up to three terms approximation for Problem \ref{TE1}.}}\label{figure5}
\end{figure}
%Figure 6 at alpha=1
\begin{figure}[H]
	\centering
	\begin{subfigure}{0.31\textwidth}
		\includegraphics[width=6.2cm]{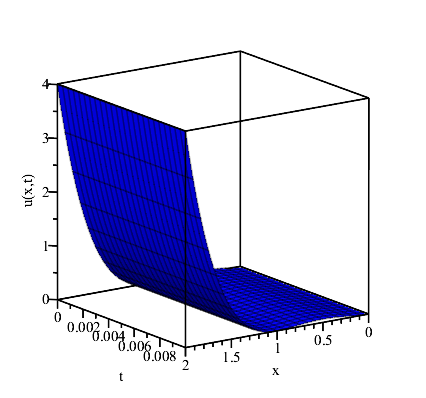}
		\caption{3D approximate  solution}
		\label{fig:noise_stabilize_u}
	\end{subfigure}\hspace{6.9em}
	%\hfill
	\begin{subfigure}{0.31\textwidth}
		\includegraphics[width=6.2cm]{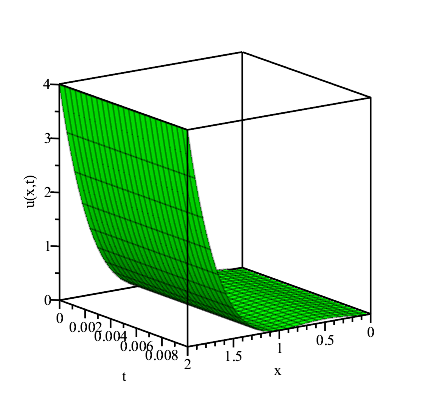}
		\caption{3D exact solution }
		\label{fig:noise_stabilize_v}
	\end{subfigure}\hspace{4.1em}
	\caption{Comparison solution plots { up to three terms approximation for Problem \ref{TE1} at $\alpha=1$.}}\label{figure6}
\end{figure}
\end{prb}

{\section{Results and Discussions}\label{sec5}
In this section, we discuss the results as obtained in Section~\ref{sec4}. We observe that the non-linearity is directly handled by using a broader concept of the Daftardar--Jafari polynomials. The numerical simulation is conducted in the following steps:
\begin{itemize}
\item[1.] The Aboodh transform (AT) is applied to the fractional derivative in the Laplace domain to obtain a simpler algebraic form of the problem.
\item[2.] The iterative procedure is considered with Daftardar--Jafari polynomials to obtain the numerical results.
\item[3.] The inverse of AT is implemented to get back to the $time$ domain from the $s$ domain.
\end{itemize}

We claim that our new scheme IATM has a higher accuracy than LADM, which can be established through numerical simulation as provided in Section~\ref{sec4}. In fact, the claimed higher accuracy is not due to the Aboodh transformation. By utilizing Daftardar--Jafari polynomials $\mathbf{J}_j$ for the non-linear term  given in \eqref{non-linear term main}, $\mathbf{J}_j$ can work more effectively as compared to Adomian polynomials $\mathbf{P}_j$. More precisely, as we have seen from the definitions of $\mathbf{J}_j$ and $\mathbf{P}_j$ given in Section~\ref{sec4}, it is obvious that the initial iterations of these polynomials are identical. However, starting from the $1^{st}$ iteration, $\mathbf{J}_{1}$ will involve more terms than $\mathbf{P}_{1}$ since more additional terms will start to emerge from $\mathbf{J}_{1}$. Due to those additional terms emerging from the Daftardar--Jafari polynomials, these terms can no doubt help achieve higher accuracy, which contributes to the difference in accuracy between the iterative schemes LADM and IATM. }

{{The beauty of fractional order cannot be underestimated, as it demonstrates a very high convergence rate towards both integer order solutions and the exact solution. We have used different fractional orders and observed that if $\alpha \in (0.5,1]$, then the solutions tend to converge to the exact solution very quickly. These results have been confirmed through subfigures \ref{figa1} and \ref{figa2}. We have also found that by using Daftardar--Jafari polynomials, the fractional order solutions converge faster compared to the fractional solutions obtained from Adomian polynomials. Thus, we conclude that IATM works accurately not only for integer order but also for fractional order solutions compared to LADM.}}

{{It is clear that Daftardar--Jafari polynomials require more computational resources compared to Adomian polynomials, but we are grateful to the STEM Lab of the Education University of Hong Kong for facilitating us with a GPU system powered by an Intel(R) Core(TM) i9-8950HK. All the experiments and plots have been conducted using Maple Version 2024 in the STEM Lab. The graphical representation confirms the effectiveness of the proposed method, where we compared the exact solution with our obtained approximate solution. We also observe that as we add more terms to the series solution $u(x,t)$, the accuracy of the methods increases gradually. It has been confirmed by the absolute error for different numbers of iterations.  }}  
\section{Conclusion and Future Work}\label{sec6}
In summary, we introduce a new method known as Iterative Aboodh transform method (IATM) in which the non-linear terms in the FPIDEs are expressed in terms of Daftardar--Jafari polynomials. The present procedure required a small number of calculations to achieve the higher accurate solutions of the targeted problems. The approximate solutions are expressed in series form of the proposed polynomials and they are illustrated with graphs and tables. It is observed that the suggested method has applicability towards both fractional and integer order problems.  { Furthermore, we find that fractional solutions gradually converge to solutions with integral order of derivative as $\alpha\to1$.} The present approach is simple to be implemented for non-linear problems, which makes it suitable for studying other non-linear FPIDEs and related systems.  { One possible extension could be the non-linear stochastic fractional integro-differential equations with suitable initial conditions. The non-linear term can be controlled by Daftardar--Jafari polynomials, and the fractional derivative can be simplified using Laplace or Abooth transformations. However, one needs to be more careful about how to control the stochastic terms in each iteration, which could be one of our possible future directions for research.}
 \subsection*{CRediT authorship contribution statement}
\textbf{Qasim Khan:} Methodology, Conceptualization, Investigation, Software, Data curation, Validation, Visualization, Resources, Formal analysis,  \& Writing – original draft. \\
\textbf{Anthony Suen:}
Supervision, Investigation, Project administration, Funding,  Writing – review \& editing.
\subsection*{Declaration of competing interest}
				The authors declare that they have no known competing financial interests or personal relationships that could have appeared to influence the work reported in this research paper.
 \subsection*{Data availability}
		No data was used for the research described in the article.
\subsection*{Funding}
This work is partially supported by Hong Kong General Research Fund (GRF) grant project number 18300821, 18300622 and 18300424, and the EdUHK Research Incentive Award project titled ``Analytic and numerical aspects of partial differential equations''.
 \subsection*{Acknowledgment }
We are very grateful to our MIT lab mates and the technicians at STEM lab of  The Education University of Hong Kong for providing all the necessary facilities for this research work.

\bibliographystyle{unsrt}
\bibliography{bib}
\end{document}